\documentclass{amsart}
\usepackage[utf8]{inputenc}

\usepackage{amsmath,amssymb,amsthm}
\usepackage{comment}
\usepackage{hyperref}
\usepackage{pgf,tikz}

\usetikzlibrary{arrows}
\usetikzlibrary{arrows.meta}

\newcommand{\RR}{\mathbb R}
\newcommand{\CC}{\mathbb C}

\renewcommand{\H}{\mathcal H}
\newcommand{\eps}{\varepsilon}

\newcommand{\abs}[1]{\left\vert #1 \right\vert}

\newcommand{\enclose}[1]{\left(#1\right)}

\newcommand{\ENCLOSE}[1]{\left\{#1\right\}}

\newcommand{\St}{\mathrm{St}}
\newcommand{\M}{\mathcal{M}}
\renewcommand{\H}{\mathcal{H}}
\newcommand{\dist}{\mathrm{dist}}
\newcommand{\co}{\mathrm{co}}
\renewcommand{\S}{\mathcal{S}}
\newcommand{\hide}[1]{}

\newtheorem{theorem}{Theorem}[section]
\newtheorem{proposition}[theorem]{Proposition}

\newtheorem{lemma}[theorem]{Lemma}

\theoremstyle{definition}

\theoremstyle{remark}
\newtheorem{remark}[theorem]{Remark}

\author[Paolini]{Emanuele Paolini} 
\address[Emanuele Paolini]{Dipartimento di Matematica, Universit\`a di Pisa \\
	Largo Bruno Pontecorvo 5 \\ I-56127, Pisa}
\email[Emanuele Paolini]{emanuele.paolini@unipi.it}

\author[Stepanov]{Eugene Stepanov}
\address[Eugene Stepanov]{%
  Universit\`a di Pisa
  \and 
  St.Petersburg Branch of the Steklov Mathematical Institute of the Russian Academy of Sciences,
	St.Petersburg, Russia
	\and
	Faculty of Mathematics, Higher School of Economics, Moscow
}
\email[Eugene Stepanov]{stepanov.eugene@gmail.com}

\thanks{
  This article was prepared within the framework of the HSE University Basic Research Program. 
  The authors acknowledge the MIUR Excellence Department Project awarded to the 
  Department of Mathematics, University of Pisa, CUP I57G22000700001.
	The first author is member of the INDAM/GNAMPA 
  and has been partially supported by Next Generation EU, Mission 4, Component 2, CUP:E53D23005860006, PRIN 2022PJ9EFL
  and by project PRA 2022 14 GeoDom (PRA 2022 - Università di Pisa) %
}
\date{\today}

\title{On the Steiner tree connecting a fractal set}

\begin{document}
\begin{abstract}
We construct an example of an infinite planar embedded self-similar binary tree  
$\Sigma$ which is the essentially unique solution to the Steiner problem of finding 
the shortest connection of a given planar self-similar fractal set $C$
of positive Hausdorff dimension.
Moreover this Steiner tree is, in fact, universal, in the sense that 
it contains, as a subtree, a Steiner tree (finite or infinite) of any given combinatorial structure.
The set $C$ can be considered the set of \emph{leaves}, or the ``boundary'',
of the tree $\Sigma$,
so that $\Sigma$ is an \emph{indecomposable} solution to the Steiner problem
with datum $C$ (i.e.\ $\Sigma\setminus C$ is connected).
\end{abstract}

\maketitle

\section{Introduction}

Let $C$ be a subset of $\RR^2$.
The Steiner problem with datum $C$ consists of finding a set $S\subset \RR^2$ 
such that $S\cup C$ is connected and $\H^1(S)$ is minimal, where $\H^1$ stands 
for the $1$-dimensional Hausdorff measure (i.e.\ the length of the set).
Usually this problem is considered for the case when $C$ is a finite 
set and, in this setting, it becomes a problem of combinatorial geometry. 
The Steiner problem goes back to Jarník and Kössler \cite{JarKos34} 
(although some of its particular cases are much more ancient, for instance 
when $C$ is the set of the three vertices of a triangle the problem was 
stated already by Fermat and solved by Torricelli).
The problem is computationally hard \cite{GarGraJoh77}.
Only few explicit examples of indecomposable solutions with an arbitrarily large number of points 
are known, namely when the set $C$ is either contained in a circle 
\cite{Rub97}
(in particular when it is the set of vertices of a regular polygon
\cite{DuHwaWen87}),
or in the case of the so called \emph{ladders} 
\cite{BurDudMai96,BurDud96,ChuGra78}.

However the problem is important for the more general setting of possibly infinite 
sets $C$. 
Moreover finding explicit solutions for an infinite $C$ may shed light on 
the solutions to finite sets with an arbitrarily large number of points.
One natural question that has arosen in this context is 
whether there exists a planar indecomposable Steiner tree $S$ connecting an uncountable 
and possibly fractal set of points $C$ 
(we say that $S$ is \emph{indecomposable} if $S\setminus C$ is 
connected).
A very similar question,
namely the existence of indecomposable Steiner sets with the topology 
of a binary tree for a finite set of points $C$ with arbitrarily 
large number of points, has been answered positively by Ivanov and 
Tuzhilin \cite{IvaTuz94} by means of an abstract construction
which gives no information on the length of the branches of the tree.
In this paper we provide a very natural construction of an infinite 
binary tree connecting a fractal set $C$ of positive
Hausdorff dimension, which is the unique solution to the Steiner problem
with datum $C$.
Moreover the tree is universal in the sense that it contains, as a subtree,
any finite or infinite Steiner tree of any given combinatorial structure
(Theorem~\ref{th:universality}).

The first explicit construction of an optimal Steiner tree connecting 
an uncountable set of vertices, has been developed
in \cite{PaoSteTep15}, were the set $C$ consists of a \emph{root} $O$ 
and uncountably many
points (\emph{leaves}) which together give a totally disconnected perfect set 
(i.e.\ every connected component is a point, but the points are not isolated).
No segment of $S$
touches the leaves, while every leaf is an accumulation point of
segments of $S$.
The infinite tree $S$ 
is composed by a \emph{trunk} of
some length $\ell$ which splits at angles of 120 degrees 
into two branches of length $\lambda_1 \ell$ 
both of which split further at the same angles into two branches of length
$\lambda_1\lambda_2\ell$ and so on (see Figure~\ref{fig:tree} right side, where the case 
$\lambda_k=\lambda$ is depicted).
\begin{figure}
  \begin{tikzpicture}
    \def\mybin#1{#1\begin{scope}
      \pgftransformxshift{1cm}\pgftransformrotate{60}\pgftransformscale{0.45}#1
      \end{scope}\begin{scope}
      \pgftransformxshift{1cm}\pgftransformrotate{-60}\pgftransformscale{0.45}#1
      \end{scope}}
    \def\myter#1{#1\begin{scope}
        \pgftransformrotate{120}#1
      \end{scope}\begin{scope}
        \pgftransformrotate{-120}#1
      \end{scope}}
    \begin{scope}
      \pgftransformscale{2}
      \node at (0.7,0) [above] {$\Sigma_u$};
      \node at (0,0) [below right] {$O$};
      \node at (1,0) [right] {$T$};
      \node at (1.65,0.5) {$A_u$};
      \myter{\mybin{\mybin{\mybin{\mybin{\mybin{\draw (0,0) -- (1,0);}}}}}}
      \begin{scope}
        \pgftransformxshift{2cm}
        \node at (0.7,0) [above] {$\Sigma$};
        \node at (0,0) [below] {$O$};
        \node at (1,0) [right] {$T$};
        \node at (1.6,0.5) {$A$};
        \mybin{\mybin{\mybin{\mybin{\mybin{\draw (0,0) -- (1,0);}}}}}
      \end{scope}
    \end{scope}
  \end{tikzpicture}
  \label{fig:tree}
  \caption{The self-similar tree with a countable 
  number of triple junctions and an uncountable number 
  of leaves which are the limit points of the triple 
  junctions. On the left the \emph{universal} tree $\Sigma_u$, 
  on the right a single branch $\Sigma$, with a fixed root in the origin.}
\end{figure}

In \cite{PaoSteTep15} it has been proven that 
the suggested tree is in fact the solution (and even the unique one) 
to the Steiner problem when the sequence of coefficients $\lambda_k$
is small and rapidly decreasing to zero (in particular summable).
In such a situation the set $C$, though being uncountable,
is not self-similar and has zero Hausdorff dimension.
The natural, original question, whether this construction is also valid 
for fractal self-similar $C$ (i.e.\ when the coefficients are constant,
i.e.\ $\lambda_j=\lambda$ for all $j$)
was left open.
In the attempt to solve it, a 
completely new and different method, based on special group valued calibrations,
to check whether the given set $S$ is a solution to the Steiner problem,
has been developed in~\cite{MarMas16a}.
In \cite{MarMas16b} the authors also show an example 
of an infinite \emph{irrigation} tree which is a solution to
a ramified transportation problem similar to the Steiner problem.
The method itself is very useful and can be applied in many situations.
However it failed in this particular problem.
A breakthrough has been achieved in the paper \cite{CheTep23}
where the authors provide a much simpler, symmetry based argument, to show that for 
sufficiently small but constant coefficients the above construction,
for a single branch with a fixed root,
gives in fact an indecomposable Steiner tree for a set $C$ 
of positive Hausdorff dimension.
In this paper we develop a method which in particular uses 
a further simplification and a generalization of 
the idea developed in \cite{CheTep23} to obtain the minimality and also uniqueness 
not only for the \emph{rooted tree} (Figure~\ref{fig:tree}, right side) 
but also for the universal fractal tree (Figure~\ref{fig:tree}, left side), 
see Theorem~\ref{th:main} and~\ref{th:universal} below.

Let us make some more precise statements.
Let $\lambda>0$ be fixed.
Consider the affine functions $f_\pm\colon \RR^2\to \RR^2$ defined by
\[
  f_\pm(x,y) = \begin{pmatrix} 1 \\ 0 \end{pmatrix} 
    + \lambda \begin{pmatrix} 
        \frac 1 2 & \frac{\mp\sqrt 3}{2} \\\\ 
        \frac{\pm\sqrt 3}{2} & \frac 1 2 
      \end{pmatrix} \begin{pmatrix} x \\ y \end{pmatrix}
\]
or, by identifying $\CC$ with $\RR^2$,
$f_\pm(z) = 1 + \lambda e^{\pm i \frac \pi 3} z$.
The function $f_+$ is the subsequent composition of: a counter-clockwise rotation of $60$ degrees,
a scaling by a factor $\lambda$, and a unit translation along the positive $x$-axis.
Let $O=(0,0)\in \RR^2$ be the origin, and let $A\subset \RR^2$ be the set of 
the limits of all the sequences obtained starting from $O$ and 
iterating $f_+$ and $f_-$.
What we obtain is a compact set $A$ which is the only non-empty compact set 
such that $A=f_+(A)\cup f_-(A)$.
If instead we start from the unit segment $[O,T]$ where $T=(1,0)$, and 
consider the union of all images of $[O,T]$ under the composition of iterations of the 
functions $f_+$ and $f_-$, we obtain the only compact set $\Sigma$ which satisfies 
$\Sigma=[O,T]\cup f_+(\Sigma)\cup f_-(\Sigma)$.
The set $\Sigma$ is depicted in Figure~\ref{fig:tree}, right side:
there is a trunk $[O,T]$ which splits into two branches $f_\pm([O,T])$, 
which in turn split into four branches $f_\pm(f_\pm([O,T]))$ and so on.
The set $A$ is composed by the \emph{leaves} of this branching tree.
If $\lambda < \frac 1 2$ (as we will assume from now on)
one can check that $A$ is 
homeomorphic to the Cantor set, it is totally disconnected, 
perfect, uncountable and it is a self-similar fractal with Hausdorff dimension
$-\log_\lambda 2<1$ (see \cite{Hut81}) and hence $\H^1(A)=0$
(the latter is in fact a necessary condition for $A$ to be the set of endpoints of topological tree
in $\RR^n$, as shown in~\cite[proposition~3.3]{PaoSteTep15}).
The set $\Sigma$ is instead a topological tree (i.e.{} connected and simply connected), 
compact, with finite positive length: $\H^1(\Sigma) = \frac{1}{1-2\lambda}$.
The set $C=\ENCLOSE{O}\cup A$ is the set of \emph{endpoints} of the tree $\Sigma$ i.e.\ 
the set of points $x\in \Sigma$ such that $\Sigma\setminus\ENCLOSE{x}$ remains connected.

The universal tree $\Sigma_u$, depicted in Figure~\ref{fig:tree}, left side,
is the union of three rotated copies of $\Sigma$: $\Sigma_u = \Sigma \cup \phi(\Sigma) \cup \phi^2(\Sigma)$
where $\phi$ is the counter-clockwise rotation of $120$ degrees centered in $O$.
Analogously we define $A_u = A \cup \phi(A) \cup \phi^2(A)$ the set of endpoints of $\Sigma_u$.

In Theorem~\ref{th:main} and \ref{th:universal},
for $\lambda$ sufficiently small,
we show that the sets 
$\Sigma$ and $\Sigma_u$ are
the shortest sets connecting $A$ and $A_u$ respectively.

%
\section{Notation and preliminaries}

For $C\subset \RR^2$ let $\St(C)$ (Steiner competitors with datum C) be the family of all sets $S$ such that 
$S\cup C$ is connected and 
let $\M(C)$ (minimal sets) be the subset of $\St(C)$ of the sets $S$ 
with minimal (possibly infinite) length $\H^1(S)$.

For a set $S\subset \RR^2$ we denote by $\bar S$ 
its closure and by $\partial S$ its topological boundary.
If $x\in \RR^2$ and $\rho>0$ then we let $B_\rho(x)\subset \RR^2$ 
be the open ball of radius 
$\rho$ centered at $x$.
When $x\in \RR^2$, and $S,S_1,S_2$ are subsets of $\RR^2$
we also write
\begin{align*}
  \dist(x,S)
    &:=\inf\ENCLOSE{\abs{x-y}:y\in S},\\
  \dist(S_1,S_2) 
    &:= \inf\ENCLOSE{\abs{x-y}\colon x\in S_1, y\in S_2}.
\end{align*}
For $P,Q\in \RR^2$ we denote by $[P,Q]$ the closed line 
segment with endpoints $P$ and $Q$, with $[P,Q)$ and $(P,Q]$ 
we denote the half-open segments, and with $(P,Q)$ we denote 
the open segment.

%
%

\subsection{Preparatory results}

In the following proposition we state some known facts about 
solutions to the Steiner problem collected from \cite{PaoSte12}
and \cite{IvaTuz94}.

\begin{proposition}[known facts about minimizers]\label{prop:PaoSte}
  Let $C\subset \RR^n$ be a compact set.
  Then $\M(C)\neq \emptyset$. If $S\in \M(C)$
  and $\H^1(S)<+\infty$
  then the following statements hold:
  \begin{itemize}
    \item[(i)] $S\cup C$ is compact, 
    and $\bar S \setminus S \subset C$;
    \item[(ii)] the closure of every connected component of $S\setminus C$     
    is a topological tree 
    (contains no subset homeomorphic to $\mathbb S^1$)
    with endpoints\footnote{%
    When $S$ is a compact topological tree, then $x\in S$ is an endpoint 
    of $S$ if and only if $S\setminus \ENCLOSE{x}$ is connected.} 
    on $C$ and with at most one endpoint 
    on each connected component of $C$;
    \item[(iii)] $S\setminus C$ is a locally finite Steiner tree in the sense that 
    for all $x\in S\setminus C$ and almost all 
    $\rho < \dist(x,C)$ the set $S\cap \bar B_\rho(x)$ 
    is the union of a finite number of line 
    segments which meet in triple points with equal 
    angles of 120 degrees and have endpoints in 
    $\partial B_\rho(x)$.

  \end{itemize}
\end{proposition}

The following lemma is also taken from~\cite[proposition~2.2 and lemma~2.6]{PaoSte12}.

\begin{lemma}\label{lm:connected}
If $S$ is connected and $\H^1(S)<+\infty$ then $\H^1(S) = \H^1(\bar S)$.
If $S$ is compact, connected, and $\H^1(S)<+\infty$ then $S$ is 
pathwise connected.
\end{lemma}

\begin{figure}
    \begin{center}
    \begin{tikzpicture}[line cap=round,line join=round,>=triangle 45,x=1.0cm,y=1.0cm]
      \draw (4.2,4.44) node {\includegraphics[scale=0.6]{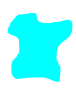}};    
      \draw(4.48,4.38) circle (1.05cm);
      \draw(4.48,4.38) circle (1.54cm);
      \draw [line width=1.6pt] (3.32,4.36)-- (3.76,4.88);
      \draw [line width=1.6pt] (3.32,4.36)-- (3.7,3.8);
      \draw [line width=1.6pt] (3.32,4.36)-- (-0.64,4.44);
      \draw [->] (4.48,4.38) -- (5.35,3.8);
      \draw [->] (4.48,4.38) -- (6.02,4.34);
      \draw (4.16,4.1) node[anchor=north west] {C};
      \draw (1.18,4.86) node[anchor=north west] {S};
      \fill (4.48,4.38) circle (1.5pt);
      \draw (4.64,4.64) node {$P$};
      \fill (3.32,4.36) circle (1.5pt);
      \draw (3.22,4.36) node[anchor=north] {$Q$};
      \fill (-0.64,4.44) circle (1.5pt);
      \draw (-0.48,4.7) node {$H$};
      \fill (2.94,4.37) circle (1.5pt);
      \draw (3.1,4.62) node {$T$};
      \draw (5.08,4.1) node {$\rho$};
      \draw (5.28,4.58) node {$r$};
      \end{tikzpicture}
    \end{center}
    \caption{The situation in Lemma~\ref{lm:angle}.}
    \label{fig:angle}
  \end{figure}

\begin{lemma}\label{lm:angle}
  Let $C\subset \RR^2$ be a compact set, $H\not \in C$ a point, 
  $C^*:=C\cup \ENCLOSE{H}$.
  Let $S\in \M(C^*)$. 
  Then there exists a point $Q\neq H$ such that $[H,Q]\subset \bar S$ 
  and either $Q\in C$ or $Q$ is a triple point of $S\setminus C^*$
  and there are no other triple points of $S\setminus C^*$ on $[H,Q]$.
  
  Moreover if $C\subset \overline {B_\rho(P)}$ for some point $P$ and radius $\rho > 0$,
  if $r\ge 2\rho/\sqrt 3$, and $r < \abs{P-H}$, then, there exists a point 
  $T\in [H,Q]\cap \partial B_r(P)$ such that 
  $\bar S\setminus B_r(P) = [H,T]$
  (see Figure~\ref{fig:angle}).
  \end{lemma}
  
\begin{proof}
  Since $H\not \in C$ there exists $\eps>0$ such that $C\cap B_\eps(H)=\emptyset$. 
  Clearly $\H^1(S\cap B_\eps(H))<+\infty$ because otherwise we could find a better competitor
  to $S$ by removing $B_\eps(H)$ from $S$ and adding the circle $\partial B_\eps(H)$ and 
  a single radius from $H$ to a point on $\partial B_\eps(H)$.
  This means that the connected components of $S\setminus C^*$
  are a finite number because each of them has length at least $\eps$ in $B_\eps(H)$.
  Let $S_1$ be a connected component of $S\setminus C^*$ such that 
  $H\in \bar S_1$.
  Recall that, by Proposition~\ref{prop:PaoSte}(iii), $S_1$ is a locally finite Steiner tree 
  composed by segments and triple points.
  We have two cases depending on whether $S_1$ contains a triple point.
  
  If $S_1$ contains triple points
  then let $Q$ be the triple point closest to $H$
  in the intrinsic (geodesic) distance on $S_1$.
  Such a point exists. 
  In fact otherwise there would exist a sequence of 
  triple points in $\bar S_1$ converging to $H$,
  because branching points can only accumulate on $C^*$ 
  by Proposition~\ref{prop:PaoSte}(ii).
  But this is not possible because in $B_\eps(H)$ there are no triple points.
  Clearly the arc from $H$ to $Q$ in $\bar S_1$ is a line segment 
  (by Proposition~\ref{prop:PaoSte}(iii)) hence $[H,Q]\subset \bar S$.

  If $S_1$ contains no triple point then $\bar S_1 =[H,Q]$ for some $Q$, 
  and, necessarily, $Q\in C$. 
  In both cases the first claim of the lemma is proven.

  For the second claim suppose that $\abs{P-H}> r$. 
  This means that the convex hull of $\ENCLOSE{H}\cup B_\rho(P)$ has an angle in 
  $H$ which is smaller than 120 degrees.
  The set $S$ must be included in this convex hull  
  and if $T$ is any triple point of $S$ then $T\in \bar B_r(P)$ because 
  otherwise at least one of the three segments meeting in $T$ would have a direction 
  which is not going to meet the ball $\bar B_\rho(P)$ and is not going towards 
  $H$. Along this segment one could hence find a further triple point $T'$ and repeat 
  the reasoning finding infinitely many triple points. 
  Since any sequence of triple points must converge to a point 
  of $C^*$, this is impossible.

  We conclude that $S\setminus \bar B_r(P)$ does not contain triple points and 
  hence is composed by a disjoint union of segments each of which has an endpoint in $H$
  and an endpoint in $\partial B_r(P)$. 
  But we can have at most one such segment because otherwise we would have, in $H$, an angle 
  smaller than 120 degrees, which is not possible.
  This means that $\bar S \setminus B_r(P)$ is contained in the segment $[H,T]$ and hence 
  is equal to a segment $[H,Q]$ with $Q \in \partial B_r(P)$.
\end{proof}

\begin{lemma}\label{lm:exists}
  Let $C$ be a compact subset of $\RR^2$ and $\ell$ be a line.
  Then $\M(\ell\cup C)\neq \emptyset$.
  If $C$ is also totally disconnected, $C\cap \ell=\emptyset$ 
  and $S\in \M(\ell\cup C)$  
  then $\bar S \supset C$.
  If moreover $\H^1(S)<+\infty$ and $\H^1(C)=0$ 
  then $\bar S\in \M(\ell\cup C)$.
\end{lemma}
\begin{proof}
  Take a compact connected segment $\sigma\subset \ell$ such that 
  the orthogonal projection of $C$ on $\ell$ is contained 
  in $\sigma$. 
  Then $\sigma\cup C$ is compact and by 
  \cite{PaoSte12} we know that 
  $\M(\sigma \cup C)\neq \emptyset$.
  For all $\Sigma \in \St(\ell\cup C)$ we have that its projection 
  $\tilde \Sigma$ 
  onto the convex hull of $\sigma\cup C$ is in $\St(\sigma\cup C)$
  and $\H^1(\tilde \Sigma)\le \H^1(\Sigma)$.
  Hence $\M(\ell\cup C)\supset \M(\sigma\cup C)\neq \emptyset$ as claimed.

  To prove the second part note that if $C$ is totally disconnected 
  then for every $x\in C$ one has $x\in \bar S$ otherwise 
  $S\cup C\cup \ell$ cannot be connected. 
  Therefore $C\subset \bar S$. 
  
  For the last claim notice that $\bar S\cap \ell$ is a finite set.
  In fact the number of connected components of $S$ arriving at $\ell$ 
  is finite since every such component must touch also points of $C$ 
  and hence has length at least $\dist(C,\ell)$.  
  Hence these components cannot be infinitely many,
  since otherwise we would have $\H^1(S)=+\infty$.
  In view of Proposition~\ref{prop:PaoSte}(ii) 
  the closure of every connected component 
  of $S\in\M(C\cup \sigma)$ 
  has at most one point on $\sigma$, showing that the set 
  $W:=\bar S\cap \sigma = \bar S \cap \ell$ 
  is finite.

  Moreover $\bar S \setminus S \subset C\cup W$ because 
  $S\cup C\cup \sigma$ is compact by Proposition~\ref{prop:PaoSte}(i).
  When $\H^1(C)=0$ this implies $\H^1(\bar S)=\H^1(S)$ and hence 
  $\bar S\in \M(\ell\cup C)$.
\end{proof}

\begin{lemma}\label{lm:base}
  Let $\ell$ be a line, $\rho>0$, 
  $r=2\rho/\sqrt 3$.
  Let $P$ be a point such that $\dist(P,\ell)>3r$
  and let $C$ be a compact subset of $\bar B_\rho(P)$.
  Let $S \in \M(\ell\cup C)$ 
  (such $S$ exists in view of Lemma~\ref{lm:exists})
  and suppose $\H^1(S)<+\infty$.
  Then there exists a point $Q\in \partial B_r(P)$ 
  and a point $H\in \ell$
  such that $\bar S\setminus B_r(P) = [H,Q]$.
  Moreover $[H,Q]$ is perpendicular to $\ell$.
\end{lemma}
\begin{proof}
  We claim that $S\setminus B_r(P)$ belongs to 
  a single connected component of $S$. 
  In fact, if not, there are two different connected components 
  $S_1$ and $S_2$ of $S$ which have some points outside of 
  $B_r(P)$. 
  Their closures must have points in $\ell$ otherwise 
  they would be contained in the convex hull of $C\subset B_r(P)$.
  By Proposition~\ref{prop:PaoSte}(ii) we know  
  that $\bar S_1$ and $\bar S_2$ each have 
  at most one endpoint on $\ell$ (hence exactly one).
  Let $H_1$ and $H_2$ be such endpoints.
  Note that $S_j\in \M(\ENCLOSE{H_j}\cup A_j)$ 
  with $A_j:=\bar S_j\cap C$, $j=1,2$, compact subsets of $C$. 
  By Lemma~\ref{lm:angle} below
  with $A_j$ and $H_j$ in place of $C$ and $H$ respectively,
  one has that
  the set $S_j\setminus B_r(P)$ is a line segment $[H_j,T_j]$
  with $T_j\in \partial B_r(P)$, $j=1,2$.
  Since $\bar S_1\cap \bar S_2\subset \ell\cup C$ 
  the line segments $[H_1,T_1]$ and $[H_2,T_2]$ may only intersect at their endpoints.
  Let $S':=(S\setminus [H_2,T_2])\cup[T_1,T_2]$ so that 
  $S'\in \St(\ell\cup C)$ and
  \[
    \H^1(S) - \H^1(S')
    = \abs{H_2 - T_2} - \abs{T_1 - T_2}
    \ge (d(P,\ell)-r) - 2r 
    > (3r-r) - 2r = 0. 
  \]
  Hence we would have $\H^1(S')<\H^1(S)$ contrary to the minimality of $S$.

  Having proven that $S\setminus B_r(P)$ belongs to a single 
  connected component of $S$, 
  denoting by $H$ its endpoint on $\ell$,
  and applying Lemma~\ref{lm:angle}, we have that 
  $S\setminus B_r(P)$ is a line segment $[H,T]$.
  Finally, if $[H,T]$ were not perpendicular to $\ell$ by moving $H$ along $\ell$ 
  we could decrease the length of $[H,T]$ contrary to the minimality 
  of $S$.
\end{proof}
  
  \begin{lemma}\label{lm:tripod}
    Let $\nu_1,\nu_2,\nu_3$ be three unit vectors in $\RR^2$
    with $\nu_1+\nu_2+\nu_3=0$. 
    Let $Y_j$ be a line perpendicular to $\nu_j$ for $j=1,2,3$.
    Let $T$ be any point in $\RR^2$ and let $H_j$ be the orthogonal 
    projection of $T$ on $Y_j$.
    Define $d_j(T) := (T-H_j)\cdot \nu_j$
    be the signed distance of $T$ from $Y_j$.
    Then $d_1(T) + d_2(T) + d_3(T)$ is constant 
    (i.e.\ independent of $T$).
  
    In particular suppose that  
    $Y_1,Y_2,Y_3$ are three lines in $\RR^2$
    forming angles of $60$ degrees so that 
    the three pairwise intersections identify
    an equilateral triangle with side length $\ell$.
    Let $T$ be any point and let $d_i(T)$ be 
    the signed distance of $T$ from $Y_i$
    with positive sign when $T$ is inside the triangle.
    Then $d_1(T) + d_2(T) + d_3(T) = \frac{\sqrt 3}{2}\ell$.
  \end{lemma}
  \begin{proof}
    If we move one of the lines parallel to itself by an amount $\delta$ 
    then the sum $d_1(T)+d_2(T)+d_3(T)$ changes by $\delta$, independent 
    of $T$.
    Therefore, without loss of generality, we may assume that the three lines 
    intersect in the origin.
    In this case $d_j(T) = T\cdot \nu_j$ and hence 
    $d_1(T)+d_2(T)+d_3(T) = T\cdot (\nu_1+\nu_2+\nu_3) = 0$
    concluding the proof of the first claim.
  
    For the second claim just notice that $d_1(T)+d_2(T)+d_3(T)$
    is constant (by the first claim) and hence is equal to the value 
    obtained when $T$ is one of vertices of the triangle, in which case 
    two distances are $0$ and the third one is equal to the height of the
    triangle, which is the claimed value.
  \end{proof}
  
\begin{lemma}\label{lm:envelope}
Let $S\in \M(\ENCLOSE{T}\cup C)$
for some $T\in \RR^2$ and a compact set $C$ contained in a horizontal strip 
$E=\ENCLOSE{(x,y)\colon \abs{y}\le \delta}$ of width $2\delta>0$.
Suppose also that there is some $T'\in S$, $T'\neq T$ 
such that $[T,T']$ is horizontal.
Then $T\in E$ and $S\subset E$.
\end{lemma}
\begin{proof}
  Suppose by contradiction that $T=(x_T,y_T)$ is outside $E$.
  For example $y_T>\delta$. 
  Then the convex hull of $\ENCLOSE{T}\cup C$ has a single 
  point, which is $T$, on the line $\ENCLOSE{y=y_T}$. 
  This is in contradiction with the fact that $T'\in S$ 
  has the same $y$-coordinate as $T$.
  This proves $T\in E$ and hence
  shows that the whole convex hull of 
  $\ENCLOSE{T}\cup C$ is contained in $E$. 
  Therefore, $S$ must be contained in $E$ as well.
\end{proof}

\begin{lemma}\label{lemma:triple}
  Let $v_1,v_2,v_3\in \RR^2$ be the three vertices of an 
  equilateral triangle of side length $\ell$. 
  Take $\rho \le \frac{\ell}{51}$.
  Let $U = B_\rho(v_1)\cup B_\rho(v_2)\cup B_\rho(v_3)$ 
  and let $A\subset U$
  be a compact set 
  such that $A_k=A\cap B_\rho(v_k) \neq \emptyset$
  for $k=1,2,3$.

  Then given any $S\in \M(A)$ there exists a point $T$ such 
  that $S=S_1\cup S_2 \cup S_3$ with $S_k\in \M(A_k\cup\ENCLOSE{T})$.
  Moreover $T\in \bar B_\rho(O)$ where $O$ is the center of the triangle.
\end{lemma}
  
\begin{proof}
Let $C$ be the convex envelope of $U$.
Since $A$ is a compact subset of $C$, given any $S\in \M(A)$ we have 
that $S$ does not touch $\partial C$ and in particular 
for all $k=1,2,3$, one has,
$S\cap \partial B_\rho(v_k) \subset \partial B_\rho(v_k) \cap C$ which 
is an arc of length $4\pi \rho/3$.
So, if we define $W_k$ to be the point on $\partial B_\rho(v_k)$ closest to $O$,
and define
\[
  S' = (S\cap U) \cup \bigcup_{k=1}^3 [O,W_k]\cup \overline{\partial B_\rho(v_k)\cap C}
\]
we clearly have $S'\in \St (A)$, 
hence
\begin{align}\label{eq:30963}
  \H^1(S) &\le \H^1(S') = \H^1(S\cap U) + 3\cdot \enclose{\frac {\ell}{\sqrt 3}-\rho + \frac 4 3 \pi \rho} \\
  &= \H^1(S\cap U) + \sqrt 3 \ell + (4\pi-3) \rho
\end{align}
We claim that $S\setminus \bar U$ contains a branching point. 
Otherwise, it would contain two disjoint paths each one connecting 
two of the three circles composing $\partial U$, hence, we would have 
\[
  \H^1(S \setminus \bar U) \ge 2(\ell-2\rho).
\]
Using~\eqref{eq:30963} we would have
\[
  2(\ell-2\rho) \le \H^1(S) - \H^1(S\cap U) \le \sqrt 3 \ell+(4\pi-3)\rho
\]
which gives is $\frac{\rho}{\ell} \ge \frac{2-\sqrt 3}{4\pi+1} > \frac{1}{51}$,
against our assumptions.

We conclude that $S'$ contains a branching point $T$. 
Let $S_1,S_2,S_3$ be the three connected components of $S\setminus \ENCLOSE{T}$.
It is clear that, up to renumbering,
one has $S_k \in \St(A_k\cup \ENCLOSE{T})$.
\end{proof}

\section{Main results}
\label{sec:main}

Let $Y:=\ENCLOSE{(x,y)\in\RR^2\colon x=0}$ be the $y$-axis,
and $Y_j := f_j(Y)$ for $j=1,2$.
Consider the point $P:=(1+\frac\lambda 2,0)$.
Clearly $\Sigma\in \St(Y\cup A)$.
We need the following lemmata.

\begin{lemma}\label{lm:existsA} 
  For every $d\in \RR$ the family $\M(\ENCLOSE{x=d}\cup A)$ is nonempty.
  If $\lambda < \frac 1 2$
  one has $\H^1(\Sigma)<+\infty$, $\H^1(A)=0$ hence $A$ is totally disconnected.
  Moreover for any $S \in \M(\ENCLOSE{x=d}\cup A)$ 
  one has that $\H^1(S)<\infty$,
  $\bar S\in \M(\ENCLOSE{x=d}\cup A)$ and
  $\bar S$ contains $A$.
\end{lemma}
\begin{proof}
  The family of minimizers is nonempty in view of Lemma~\ref{lm:exists}.

  We claim that since $\lambda < \frac 1 2$ then $\H^1(A)=0$. 
  Since $A\subset \Sigma$ and $\H^1(\Sigma)=\frac{1}{1-2\lambda}$
  is finite if $\lambda<\frac 1 2$ then $\H^1(A)<+\infty$.
  By the self-similarity property of $A$ we have that 
  \[
  \H^1(A) = 2 \lambda \H^1(A).
  \]
  If $0<\H^1(A)<+\infty$ we would have $\lambda = \frac 1 2$.
  Hence, if $\lambda < \frac 1 2$ we conclude that $\H^1(A)=0$.
  This also means that $A$ is totally disconnected because each connected 
  component of $A$ is pathwise connected (by Lemma~\ref{lm:connected}) 
  and hence if it were not a single point it would have positive length.
    
  Now the conclusion follows from Lemma~\ref{lm:exists} with $C:=A$ 
  minding that 
  \[
    \H^1(S)\le \H^1(\Sigma) = \frac{1}{1-2\lambda}<+\infty.
  \]
\end{proof}

\begin{lemma}\label{lm:01}
  Let $T_0=(1,0)$ and $A$ be defined as above,
  $\lambda \le \frac 1 8$,
  $\rho= \frac{\lambda}{1-\lambda}$,
  $r=\frac{2\rho}{\sqrt 3}$
  Then $A\subset B_\rho(T_0)$ and 
  given any $d\le \frac 1 2$ and any
  $S\in \M(\ENCLOSE{x=d}\cup A)$ 
  there exists $T\in \partial B_r(T_0)$ such that
  we have that $S\setminus B_r(T_0)$ is a line segment
  $[H,T]$ perpendicular to $\ENCLOSE{x=d}$ 
  with $H\in \ENCLOSE{x=d}$.
  In particular $\bar S$ is connected.
\end{lemma}
\begin{proof}
  The claim $A\subset B_\rho(T_0)$ can be easily checked
  by noticing that any point of $A$ has distance
  from $T_0$ smaller than 
  \[
      \lambda + \lambda^2 + \lambda^3 + \dots = \rho.
  \]
  Lemma \ref{lm:base} with $C:=A$, $P:=T_0$, $\ell:=\ENCLOSE{x=d}$ 
  can be applied since the assumption $\lambda \le 1 / 8$ assures 
  that $3r < 1/2 \le 1-d = \dist(T_0,\ENCLOSE{x=d})$.
  So the point $T$ exists with the desired properties and 
  $\bar S\cap \ENCLOSE{x=d} =\ENCLOSE{H}$.
  Since, by Lemma~\ref{lm:existsA}, $\bar S \supset A$
  it follows that $\bar S \supset \ENCLOSE{H}\cup A$ which is totally disconnected 
  and hence $\bar S$ must be connected.
\end{proof}

\begin{lemma}\label{lm:precedente1}
One has
\begin{align*}
    \dist(Y,f_j(A)) 
    &=\dist(Y,A)
    \ge
    1+ \frac{\lambda} 2 
    - \frac{\lambda^2(1+2\lambda)}{2(1-\lambda^2)},
    \qquad j=1,2,
    \\
    \dist(f_1(A),f_2(A)) & 
    \ge 
    \sqrt 3 \lambda - \sqrt 3 \frac{\lambda^3}{1-\lambda}.
\end{align*}
\end{lemma}
\begin{proof}
If we follow the tree $\Sigma$ starting from the root 
at the origin and trying to get as close as possible 
to the line $Y$ we first must follow the trunk of length 
$1$ then at the first branching nothing changes if we
turn left or right and then we go further along 
the oblique branch of length $\lambda$ so that the distance 
from $Y$ increases by $\lambda /2$. 
After that the path leading towards $Y$ is made 
of alternating branches with horizontal and oblique directions.
The computation gives
  \begin{align*}
  \dist(Y,f_j(A)) &\ge 1+ \frac \lambda 2 
     - \enclose{\frac{\lambda^2}{2} 
      + \lambda^3 
      + \frac{\lambda^4 }{2}
      + \lambda^5
      + \frac{\lambda^6}{4} + \dots}\\
      &= 1 + \frac \lambda 2 
      - \frac{1}{2}\frac{\lambda^2}{1-\lambda^2}
      - \frac{\lambda^3}{1-\lambda^2}\\
      &= 1+ \frac{\lambda} 2 
      - \frac{\lambda^2(1+2\lambda)}{2(1-\lambda^2)}.
  \end{align*}
  Analogously,  
  \begin{align*}
    \dist(f_1(A),f_2(A))
    &\ge 2\, \dist(f_1(A),\ENCLOSE{y=0}) \\
    &\ge 2\enclose{\frac{\sqrt 3}{2} \lambda 
      -\frac{\sqrt 3}{2}\enclose{\lambda^3 + \lambda^4 + \dots}
     }
    = \sqrt 3 \lambda -\sqrt 3 \frac{\lambda^3}{1-\lambda},
  \end{align*}
  as claimed.
\end{proof}

\begin{remark}
  Notice that if the right hand side of the inequalities 
  in the claim of Lemma~\ref{lm:precedente1}
  are positive then in fact those inequalities 
  become equalities as it is easy to see from the proof.
  This happens for $\lambda$ sufficiently small, 
  for instance when $\lambda < \frac 1 2$.
\end{remark}

\begin{figure}
  \begin{center}
  \begin{tikzpicture}[line cap=round,line join=round,>=triangle 45,x=4.0cm,y=4.0cm]
    \clip(-0.2,-0.7) rectangle (1.84,0.92);
    \draw [color=yellow,line width=3](0,0.71334)--(1.24,0)--(0,-0.71334)--cycle;
    \draw (1.0,0.57) node[anchor=north west] {\includegraphics[width=1cm]{blob1.png}};
    \draw (1.19,0.5) node[anchor=north west] {$f_1(A)$};
    \draw (1.05,-0.12) node[anchor=north west] {\includegraphics[width=1cm]{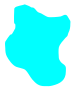}};
    \draw (1.22,-0.28) node[anchor=north west] {$f_2(A)$};
    \draw (1.2,0.35)-- (1,0);
    \draw (1,0)-- (1.2,-0.35);
    \draw [line width=1.2pt] (0,0)-- (1,0);
    \draw (0,0.58)-- (0,-0.58);
    \draw [domain=-0.2:1.84] plot(\x,{(-0.58--0.58*\x)/1});
    \draw [domain=-0.2:1.84] plot(\x,{(-0.58--0.58*\x)/-1});
    \draw (0,0)-- (1,0);
    \draw [line width=1.2pt] (1.2,0.35)-- (1,0);
    \draw [line width=1.2pt](1,0)-- (1.2,-0.35);
    \draw (0.1,-0.7) -- (0.1,0.92);
    \draw [line width=1.2pt,color=blue] (1.15,0.37)-- (1,0.1);
    \draw [line width=1.2pt,color=blue] (1,0.1)-- (1.21,-0.27);
    \draw [line width=1.2pt,color=blue] (1,0.1)-- (0.1,0.1);
    \draw [dash pattern=on 1pt off 1pt,domain=-0.2:1.84] plot(\x,{(-0.68--0.58*\x)/-1});
    \draw [dash pattern=on 1pt off 1pt,domain=-0.2:1.84] plot(\x,{(-0.47--0.58*\x)/1});
    \draw [dash pattern=on 2pt off 1pt,domain=-0.2:1.84] plot(\x,{(-0.71334--0.58*\x)/-1});
    \draw [dash pattern=on 2pt off 1pt,domain=-0.2:1.84] plot(\x,{(-0.71334--0.58*\x)/1});
    \draw (0.25,0.1)-- (0.25,0);
    \draw (0.1,0.22)-- (0,0.22);
    \draw [{Latex[length=1mm]}-{Latex[length=1mm]}] (0,0.22) -- (0.1,0.22);
    \draw [{Latex[length=1mm]}-{Latex[length=1mm]}] (0.25,0.1) -- (0.25,0);
    \draw [{Latex[length=1mm]}-{Latex[length=1mm]}] (0.67,0.29) -- (0.63,0.22);
    \draw [{Latex[length=1mm]}-{Latex[length=1mm]}] (0.65,-0.1) -- (0.7,-0.18);
    \draw [{Latex[length=1mm]}-{Latex[length=1mm]}] (1.05,0.2) -- (1.1,0.17);
    \draw [{Latex[length=1mm]}-{Latex[length=1mm]}] (1.09,-0.07) -- (1.05,-0.09);
    \draw [domain=-0.2:1.84] plot(\x,{(-0-0*\x)/1});
    \draw (0,-0.7) -- (0,0.92);
    \begin{scriptsize}
    \fill (0,0) circle (1.5pt);
    \draw (0,0) node[anchor=north east] {$O$};
    \fill (1,0) circle (1.5pt);
    \draw (0.995,-0.00) node[anchor=north] {$T_0$};
    \fill (1.23,0) circle (1.5pt);
    \draw (1.24,0.0) node[anchor=south] {$P$};
    \draw (0,0.4) node[anchor=east] {$Y$};
    \draw (0.21,0.4) node[anchor=east] {$Y'$};
    \draw (-0.17,-0.63) node {$Y_2$};
    \draw (-0.17,0.63) node {$Y_1$};
    \fill [color=blue] (0.1,0.1) circle (1.5pt);
    \draw (0.1,0.1) node[anchor=south west] {$H$};
    \fill [color=blue] (1,0.1) circle (1.5pt);
    \draw (1.01,0.12) node[anchor=south] {$T$};
    \draw (-0.16,0.715) node {$Y'_1$};
    \draw (-0.17,-0.51) node {$Y'_2$};
    \draw (0.45,0.50) node {$Y''_1$};
    \draw (0.45,-0.51) node {$Y''_2$};
    \draw (0.36,0.05) node {$\delta=\delta_0$};
    \draw (0.12,0.26) node {$d=d_0$};
    \draw (0.61,0.28) node {$d_1$};
    \draw (0.71,-0.12) node {$d_2$};
  \end{scriptsize}
  \begin{tiny}
    \draw (1.09,0.21) node {$\delta_{\!1}$};
    \draw (1.1,-0.12) node {$\delta_{\!2}$};
  \end{tiny}
    \end{tikzpicture}
    \caption{Constructions and notation in the proof of 
    Lemma~\ref{lm:branching} and Theorem~\ref{th:main}.
    The triangle $\Delta$ from the proof of Lemma~\ref{lm:branching} is highlighted.}
  \end{center}
\end{figure}

\begin{lemma}\label{lm:branching}
  Let $Y':=\ENCLOSE{x=d}$ for some $d\le\frac 1 2$,
  be a line parallel to $Y=\ENCLOSE{x=0}$.
  If $S\in \M(Y'\cup A)$
  and $\lambda < \frac 1 {25}$ then
  there is a branching point $T\in S$ 
  such that $S\setminus \ENCLOSE{T}$ is composed of 
  three connected sets, the closures of which are
   $[H,T]$, $S_1$ and $S_2$ respectively,
  with $H\in Y'$ and $[H,T]$ perpendicular to $Y'$.  
  Moreover, 
  \begin{equation}\label{eq:474947}
    \dist(H,\ENCLOSE{y=0})\le 2\lambda\rho
    \qquad\text{and}\qquad 
    T\in B_{2\lambda\rho}(T_0), 
  \end{equation}
  where $\rho:=\frac{\lambda}{1-\lambda}$ and $T_0:=(1,0)$.
  Finally if $Y_j'$ is the line parallel to $Y_j:=f_j(Y)$ 
  passing through $T$ 
  then $S_j \in \M(Y_j'\cup f_j(A))$ 
  and $Y_j'\subset f_j(\ENCLOSE{x\le\frac 1 2})$,
  for $j=1,2$.
\end{lemma}
\begin{proof}
  \emph{Step 1.}
  We first claim that $\bar S$ cannot contain 
  two compact connected sets 
  $\sigma$ and $\gamma$
  with $\H^1(\sigma\cap \gamma)=0$
  such that $\sigma$ touches both $Y'$ and $f_1(A)$ 
  while $\gamma$ touches both $f_1(A)$ and $f_2(A)$.
  In fact, otherwise we would have
  \[
    \H^1(\bar S)
    \ge \H^1(\sigma) + \H^1(\gamma)
    \ge \dist(Y', f_1(A)) + \dist(f_1(A), f_2(A)),
  \]
  but $\dist(Y', f_1(A)) \ge \dist(Y,f_1(A)) - d$,
  hence,
  by Lemma~\ref{lm:precedente1},
  we would get
  \begin{equation}
  \label{eq:43747}
  \begin{aligned}
    \H^1(\bar S)
    \ge 
    1 - d + \frac{\lambda} 2 
      - \frac{\lambda^2(1+2\lambda)}{2(1-\lambda^2)}
     +
     \sqrt 3 \lambda - \sqrt 3 \frac{\lambda^3}{1-\lambda}.
  \end{aligned}
  \end{equation}
  Let $O'=(d,0)$ and $\Sigma':=[O',T_0] \cup f_1(\Sigma) \cup f_2(\Sigma)$ 
  be the tree obtained by adding or removing a segment of length 
  $\abs{d}$ from $\Sigma$ so that $\Sigma'\in \St(Y'\cup A)$.
  We have  
  \begin{equation}\label{eq:44321}
    \H^1(\Sigma')
    = \H^1(\Sigma) - d 
    = 1 - d + 2 \lambda + 4 \lambda^2 + \dots 
    = \frac{1}{1-2\lambda}-d
  \end{equation}
  and one can check that the rhs of \eqref{eq:43747} is
  strictly greater than the rhs of \eqref{eq:44321}
  for $\lambda < \frac 1 {25}$.
  Hence $\H^1(S) =\H^1(\bar S) > \H^1(\Sigma')$ contrary 
  to the optimality of $S$.
  
  \emph{Step 2.} 
  By Proposition~\ref{prop:PaoSte}(ii) we know that $S$ touches $Y'$ in a single point $H$.
  Recall that $\bar S$ is connected in view of Lemma~\ref{lm:01}.
  Consider an injective arc $\theta\colon[0,1]\to \bar S$ 
  such that $\theta(0)=H$, $\theta([0,1))\cap A=\emptyset$ 
  and $\theta(1)\in A$.  
  Let 
  \[
    t:= \sup\ENCLOSE{s\in[0,1]\colon \bar S\setminus \theta([0,s])
  \text{ is connected}}
  \] 
  and let $T:=\theta(t)$, $S_0:=\theta([0,t])$. 
  Observe that $(\bar S\setminus S_0)\cup \ENCLOSE{T}$ 
  is compact and arcwise connected.
  In fact for every pair of points $P,Q\in (\bar S\setminus S_0)\cup \ENCLOSE{T}$
  there exists a unique injective arc $\Gamma$ in $\bar S$ 
  connecting $P$ and $Q$ (this follows by Proposition~\ref{prop:PaoSte}(ii)
  since $\bar S$ is connected and having finite length 
  is an arcwise connected topological tree).
  It is enough to note now that 
  $\Gamma\cap \theta([0,s])=\emptyset$
  for all $s<t$ and hence 
  $\Gamma\subset \bar S\setminus \theta([0,t))
  = (S\setminus S_0)\cup\ENCLOSE{T}$.
  
  \emph{Step 3.} 
  We are going to define $S_1$ and $S_2$ closed connected 
  subsets of $\bar S$ such that 
  $\bar S = S_0\cup S_1 \cup S_2$ 
  and $S_j$, $j=0,1,2$, have $T$ 
  as the only common point, and moreover,
  $S_j$ contains $f_j(A)$ for $j=1,2$.
  To this aim we are going to consider two cases 
  which will be excluded.

  \emph{Case 1.} 
  Suppose $T$ is a point of $A$.
  Without loss of generality suppose $T\in f_1(A)$.
  The set $(\bar S\setminus S_0)\cup\ENCLOSE{T}$
  is connected and contains points of both $f_1(A)$
  and $f_2(A)$ 
  because $(S_0\setminus\ENCLOSE{T})\cap A = \emptyset$ 
  by construction hence there exists an arc $\gamma$ 
  in $(\bar S\setminus S_0)\cup\ENCLOSE{T}$.
  The claim of Step~1 with $\sigma:=S_0$ implies that this cannot happen.
  
  \emph{Case 2.} Suppose $T\not \in A$.
  In this case $S$ is a finite tree in a neighbourhood of $T$
  by Proposition~\ref{prop:PaoSte}(iii) 
  and therefore $T$ is a triple point of $S$. 
  Hence $\bar S\setminus S_0$ has two connected components 
  $S_1'$ and $S_2'$ (recall that $\bar S$ contains no loops).
  Let $S_j:=S_j'\cup \ENCLOSE{T}$, $j=1,2$.
  Each $S_j$ must contain at least one point of $A$ because otherwise 
  $\bar S\setminus S_j'\in \St(Y\cup A)$ will be 
  strictly shorter than $S\in \M(Y\cup A)$.
  Moreover each point of $A$ is contained in either $S_1$ or $S_2$.
  Without loss of generality suppose that $S_1$ 
  contains at least one point of $f_1(A)$.
  
  \emph{Case 2a.} 
  Suppose that $S_1$ contains also points of $f_2(A)$.
  Then we can apply the claim of Step~1 with 
  $\sigma:= S_0 \cup S_2$ and 
  $\gamma:= S_1$ and exclude that this case can happen.
  
  \emph{Case 2b.} If $S_2$ contains points of both 
  $f_1(A)$ and $f_2(A)$ we proceed 
  as in Case 2a with $S_1$ and $S_2$ interchanged.
  
  The only remaining possibility is that $S_1$ 
  only touches points of $f_1(A)$ 
  and $S_2$ only touches points of $f_2(A)$. 
  Hence $S_1\supset f_1(A)$ and $S_2\supset f_2(A)$ since $A\subset \bar S$.
  Therefore, $S_0\in \M(Y'\cup \ENCLOSE{T})$  
  and $S_j \in \M(\ENCLOSE{T} \cup f_j(A))$ for $j=1,2$:
  otherwise, by substituting $S_0$ with an element of $\M(Y'\cup \ENCLOSE{T})$
  and $S_j$ with an element of $\M(\ENCLOSE{T}\cup f_j(A))$,
  we could construct a better competitor 
  to $S\in \St(Y'\cup A)$.
  
  \emph{Step 4.} 
  Clearly $S_0\in \M(Y'\cup \ENCLOSE{T})$ implies that $S_0$ is the straight 
  line segment $[H,T]$ perpendicular to $Y'$.
  
  \emph{Step 5.}
  By Lemma~\ref{lm:01} we have $A\subset B_\rho(T_0)$,
  with $\rho=\frac{\lambda}{1-\lambda}$
  hence $f_j(A) \subset B_{\lambda\rho}(f_j(T_0))$
  for $j=1,2$.
  As noticed in Step~3 we know that $S$ is a regular tripod in a small 
  neighbourhood of $T$ composed by three line segments forming 
  equal angles of 120 degrees. Since $S_0$ is perpendicular to $Y$ 
  it follows that $S_j$ contains a small segment perpendicular to $Y_j=f_j(Y)$.
  Since $S_j\in \M(\ENCLOSE{T}\cup f_j(A))$ for $j=1,2$,
  by Lemma~\ref{lm:envelope} (applied to $S_j$ in place of $S$, 
  $f(A_j)$ in place of $C$ and coordinate $y$ along the axis $Y_j$ 
  and $x$ perpendicular to $Y_j$),
  we obtain that $S_j$ is contained in the strip perpendicular 
  to $Y_j$, centered in $T_0$ and containing $f_j(A)$. 
  Since $f_j(A)\subset B_{\lambda\rho}(f_j(T_0))$
  such a strip has width $2\lambda\rho$
  (notice that $f_j(T_0)$ lies on the line passing through $T_0$ and perpendicular 
  to $Y_j$).
  
  The intersection between the two strips for $j=1,2$ 
  is the union of two equilateral 
  triangles each with height $2\lambda\rho$. 
  Hence $T$ is contained in the ball
  $B_{2\lambda\rho}(T_0)$, 
  and thus, in particular, 
  \[
    \dist(H,\ENCLOSE{y=0})= \dist(T,\ENCLOSE{y=0})\le \abs{TT_0} \le 2\lambda\rho,
  \]
  proving~\eqref{eq:474947}.
  Moreover, also the distance of the line $Y_j'$ from $T_0$ is less than 
  $2\lambda\rho$ hence $Y_j'\subset f_j(\ENCLOSE{x<2\rho})
  \subset f_j(\ENCLOSE{x<\frac 1 2})$, for $j=1,2$, since $\lambda < \frac 1 5$.
  
  \emph{Step 6.}
  We claim that each $S_j$, $j=1,2$, has no branching points inside the triangle 
  delimited by $Y$, $Y_1$ and $Y_2$ 
  (to avoid confusion, notice that $T$ is not a branching point of $S_j$).
 
  If $T$ is itself outside of this triangle there is nothing to prove
  because $S_j\subset \overline{\co}(\ENCLOSE{T}\cup A_j)$ 
  is also outside of the triangle.
  Otherwise, 
  since $f_j(A)\subset B_{\lambda\rho}(f_j(P))$,
  then, by Lemma~\ref{lm:base} all branching points of $S_j$ 
  are inside $B_{\frac{2\lambda \rho}{\sqrt 3}}(f_j(T_0))$ while 
  $d(f_j(T_0),Y_j)=\lambda d(T_0,Y)
  = \lambda > \frac{2\lambda \rho}{\sqrt 3}$
  since $\lambda < \frac{3}{2\sqrt 3+3}\approx 0.46$.

  \emph{Step 7.}
  Let $Y_j''$ be the line parallel to $Y_j$ and passing through $P$.
  Denote with $\Delta$ the triangle delimited by $Y$, $Y_1''$ and $Y_2''$,
  and notice that $T$ is inside $\Delta$ because $T\in B_{\lambda^2}(T_0)$
  (as proven in Step~5)
  while $\dist(T_0, Y_j'') = \frac{\lambda}{4}>2\lambda\rho$ for $j=1,2$.
  Let $T_j''$ be the point on $Y_j''$ such that $[T,T_j'']$ is perpendicular
  to $Y_j$. 
  Recall that $S_j$ in a neighbourhood of $T$ is a segment perpendicular 
  to $Y_j$ because $S$, in a neighbourhood of $T$, 
  is a regular tripod with equal angles of 120 degrees 
  and $S_0$ is perpendicular to $Y$.
  We claim that $[T,T_j'']$ is contained in $S_j$: otherwise there would 
  be a branching point of $S_j$ on such segment, which is not the case because,
  as already stated, all branching points of $S_j$ are inside 
  the ball $B_{\frac{2\lambda \rho}{\sqrt 3}}(f_j(T_0))$, while
  \begin{align*}
    \dist(f_j(A), Y_j'') 
    &\ge \dist(f_j(A), Y_j) - \dist(Y_j, Y_j'')
    \ge \lambda \dist(A,Y) - \frac{\lambda}{4}  \\
    & \ge \lambda \enclose{1 + \frac \lambda 2 + \frac{\lambda^2}{1-\lambda}} - \frac \lambda 4 \\ 
    & \ge \frac 3 4 \lambda
    > \frac 2 {\sqrt 3}\lambda\rho
    \qquad \text{for $\lambda < \frac{9}{8\sqrt 3+9}\approx 0.39$}. \\
  \end{align*}
  Hence $\tilde S_j:= S_j\setminus [T,T_j'']$ is connected and 
  $\H^1(\tilde S_j)=\H^1(S_j) - \dist(T,Y_j'')$.

  \emph{Step 8: conclusion.}
  We now show the last claim, $S_j \in \M(Y_j'\cup f_j(A))$, $j=1,2$.
  To this aim take any $S_j''\in \M(Y_j''\cup f_j (A))$ and let $H_j''$ be the 
  endpoint of $S_j''$ on the line $Y_j''$ ($S_j''$ touches $Y_j''$ in a single point 
  in view of Lemma~\ref{lm:angle}). 
  Note that $H_j''$ is on the boundary of the triangle $\Delta$ 
  in view of Lemma~\ref{lm:envelope} 
  because the orthogonal 
  projection of $f_j(A)$ onto $Y_j''$ belongs to the respective side of this triangle.
  Define $S_j':= S_j''\cup [H_j'',H_j']$ where $H_j'$ is the projection of $H_j''$ 
  onto the line $Y_j'$, $j=1,2$ and let $\Gamma\in \M(Y\cup \ENCLOSE{H_1'',H_2''})$
  be the regular tripod with endpoints in $H_1''$, $H_2''$ and a third point on 
  the side of $\Delta$ contained in $Y$.
  Clearly $S_j'\in \M(Y_j'\cup f_j(A))$ and $\H^1(S_j')=\H^1(S_j'') + \abs{H_j H_j''}$
  while $\H^1(\Gamma) = 1+\frac \lambda 2$ in view of Lemma~\ref{lm:tripod}.

  Therefore, if we consider $\tilde S := \Gamma \cup S_1''\cup S_1''$ 
  we have $\tilde S\in \St(Y\cup A)$ hence $\H^1(\tilde S)\ge \H^1(S)$.
  On the other hand, we have 
  \begin{align*}
    \H^1(\tilde S) 
    &= \H^1(\Gamma) + \H^1(S_1'') + \H^1(S_2'') \\
    &\le 1 + \frac{\lambda}{2} 
      + \H^1(S_1') - \abs{H_1' H_1''}
      + \H^1(S_2') - \abs{H_2' H_2''} \\
      &\le 1 + \frac \lambda 2 + \H^1(S_1) - \abs{H_1' H_1''} + \H^1(S_2) - \abs{H_2' H_2''}\\
      &= 1  + \frac \lambda 2 + \H^1(S_1) - \abs{T T_1''} + \H^1(S_2) - \abs{T T_2''}\\
    \end{align*}
    but $\H^1(S_0)+\abs{T T_1''}+\abs{T T_2''} = 1+\frac \lambda 2$
    by Lemma~\ref{lm:tripod}
    and hence
    \[
      \H^1(\tilde S) 
      \le \H^1(S_0) + \H^1(S_1) + \H^1(S_2)
      = \H^1(S).
    \]
  Hence, all the above inequalities are, in fact, equalities,
  which means in particular that $\H^1(S_j')=\H^1(S_j)$
  and thus $S_j \in \M(Y_j'\cup f_j(A))$, $j=1,2$ as claimed.
\end{proof}

Finally we are in the position to prove the minimality of 
the branching tree $\Sigma$ with a fixed root, or, more generally,
with the root moving along the line $Y$.

\begin{theorem}\label{th:main}
Suppose $\lambda < \frac 1{25}$.
Then one has $\Sigma \in \M(\ENCLOSE{0}\cup A) = \M(Y\cup A)$.
Moreover $S\in \M(Y\cup A)$ if and only if 
$\Sigma\setminus (\ENCLOSE{0}\cup A)\subset S\subset \Sigma$,
which in particular implies $\bar S=\Sigma$.
\end{theorem}
\begin{proof}
By Lemma~\ref{lm:existsA} we know that $\M(\ENCLOSE{x=d}\cup A)$
is nonempty for every $d\in \RR$ and
that for any $S\in \M(\ENCLOSE{x=d}\cup A)$ 
one has that $\bar S\in \M(\ENCLOSE{x=d}\cup A)$ and $\bar S\supset A$.
Moreover, by Lemma~\ref{lm:01},
$\bar S$ is connected.
Thus we may consider an arbitrary closed $S\in \M(\ENCLOSE{x=d}\cup A)$.
If $d<\frac 1 2$, 
by Lemma~\ref{lm:branching}, 
we are able to define $S_1$ and $S_2$, $H$, $T$
such that 
$S=[H, T] \cup S_1 \cup S_2$, $[H,T]$ is perpendicular to $Y$.
If we define the sets $b_1(S)$ and $b_2(S)$ 
(the ``\emph{branches}'' of $S$) 
by $b_j(S) := f_j^{-1}(S_j)$, $j=1,2$.
We notice that $b_j(S)\in \M(\ENCLOSE{x=d_j}\cup A)$ for some $d_j<\frac 1 2$, 
as claimed in the Lemma~\ref{lm:branching}, $j=1,2$.
Let $d(S)$ and $\delta(S)$ be the two coordinates of the point $H$ 
so that $d(S)$ is the distance of $H$ from the line $Y$ 
and $\delta(S)$ is the distance of the point $H$
(or, which is the same, the point $T$) 
from the line $X:=\ENCLOSE{y=0}$.
Define $\S^0:=\ENCLOSE{\tilde S}$ where $\tilde S\in \M(\ENCLOSE{x=0}\cup A)$ 
is fixed, compact set,
and define inductively $\S^k$ as the family of the rescaled branches of $\tilde S$ 
at level $k$:
\[
  \S^{k+1} := \bigcup_{S\in \S^k}\ENCLOSE{b_1(S),b_2(S)}.
\]
Notice that Lemma~\ref{lm:branching} can be applied to $\tilde S$ and 
hence inductively on the two rescaled branches 
$b_1(S)$ and $b_2(S)$ for every $S\in \S^k$
since the properties 
$S \in \M(\ENCLOSE{x=d(S)}\cup A)$ and $d(S)<\frac 1 2$
are preserved by the operators $b_1$ and $b_2$
as stated in Lemma~\ref{lm:branching} itself.

We claim that for all $S\in \S^k$, one has 
\[
  \delta(S)\le 2\lambda \max\ENCLOSE{\delta(b_1(S)),
    \delta(b_2(S))}.
\]
We know that $S\in \S^k$ touches the vertical line $\ENCLOSE{x=d(S)}$ 
in a single point $H$.
By Lemma~\ref{lm:branching} we know that $S$ is composed 
by a segment $[H, T]$ and two trees $S_1$, $S_2$
with $S_j\in \M(Y_j'\cup f_j(A))$, $j=1,2$
where $Y_j'$ is the line passing through $T$ and parallel to $Y_j$.
Let $\delta_j$ be the distance of $T$ from $f_j(X)$, $j=1,2$.
By Lemma~\ref{lm:branching} one has $\delta_j \le 2\lambda\rho$.
By Lemma~\ref{lm:tripod} applied to the three lines $X$, $f_1(X)$ and $f_2(X)$
(all lines passing through $T_0$) we know that the sum of the signed distances of $T$ from the three lines 
passing through $T_0$ is equal to $0$.
The absolute values of these three distances
are, respectively, $\delta(S)$, $\delta_1$ and $\delta_2$,
so that
\begin{align*}
  \delta(S) &\le \delta_1 + \delta_2 
    \le 2\max\ENCLOSE{\delta_1,\delta_2}\\
    &= 2\max\ENCLOSE{\lambda\delta(b_1(S)),\lambda\delta(b_2(S))}
  \end{align*}
because $b_j(S)=f_j^{-1}(S_j)$,
showing the claim.

The proven claim now gives 
\[
  \Delta_k 
  := \max\ENCLOSE{\delta(S)\colon S\in \S^k}
   \le 2\lambda \max\ENCLOSE{\delta(S)\colon S\in \S^{k+1}}
   = 2 \lambda \Delta_{k+1}.
\]
We have $\Delta_0 = 0$ since otherwise for 
$\Delta_0>0$ we would have 
$\Delta_k\to +\infty$ (because $2\lambda <1$), while we know,
by Lemma~\ref{lm:branching}, that $\Delta_k \le 2\lambda\rho$
for all $k$.
But $\Delta_0 =0$ implies $\Delta_k=0$ for all $k$.
In particular for $d=0$ this implies $\tilde S=\Sigma$ 
(recall that $\tilde S$ has been choosen to be closed).

Note that we have proven that for any arbitrary $S\in \M(Y\cup A)$ 
one has $\bar S= \Sigma$ and in particular $S\in \M(\ENCLOSE 0 \cup A)$. 
But we also have that $\bar S \setminus S \subset \ENCLOSE{0}\cup A$
by Proposition~\ref{prop:PaoSte}(i),
that is $\Sigma\setminus S \subset \ENCLOSE 0 \cup A$.
This implies that $S\supset \Sigma\setminus(\ENCLOSE 0 \cup A)$.

It remains to prove that $\M(\ENCLOSE 0 \cup A) = \M(Y\cup A)$.
To this aim note that we have already proven 
that $\M(Y\cup A) \subset \M(\ENCLOSE 0 \cup A)$
but for every $S\in \M(\ENCLOSE 0 \cup A)$ one necessarily has 
$S\in \St(Y\cup A)$ and $\H^1(S)=\H^1(\Sigma)$ while
$\Sigma \in \M(Y\cup A)$ which implies the claim.

\end{proof}

We can now prove the main result for the universal tree $\Sigma_u$.

Let $A$, $\Sigma$, $A_u$, $\Sigma_u$, be the sets defined in the introduction
(see Figure~\ref{fig:tree})
$A_u = A \cup \phi(A)\cup \phi^2(A)$ and 
$\Sigma_u = \Sigma \cup \phi(\Sigma)\cup \phi^2(\Sigma)$
where $\phi$ is the counter-clockwise rotation of angle of 120 degrees.

\begin{theorem}[universal minimizer]\label{th:universal}
  Suppose $\lambda < \frac 1{29}$.
  Then, $S \in \M(A_u)$ if and only if $S\cup A_u = \Sigma_u$.
\end{theorem}
\begin{proof}
  Let $v_1 = (1,0)\in \RR^2$, $v_2=\phi(v_1)$, $v_3=\phi^2(v_1)$ 
  be the vertices of an equilateral 
  triangle of side length $\ell = \sqrt 3$ centered in the origin.
  Let $S\in \M(A_u)$.
  Let $\rho = \frac{\lambda}{1-\lambda}$, so that $A\subset B_\rho(v_1)$ as 
  stated in Lemma~\ref{lm:01}.
  By Lemma~\ref{lemma:triple} we know that there exists a point 
  $T\in B_\rho(0)$ such that $S=S_1\cup S_2\cup S_3$ with $S_k\in \M(A_k\cup \ENCLOSE{T})$.
  Consider the triangle $\Delta$ 
  with vertices $\enclose{\frac 1 2,\pm\frac{\sqrt 3}{2}}$ and $\enclose{-1,0}$.
  Let $Y$ be the line with equation $x=\frac 1 2$.
  Clearly $\phi(\Delta)=\Delta$ and $\Delta$ has a side contained in the line $Y$. 
  By Lemma~\ref{lm:angle} we know that $S_1\setminus B_r(V_1)$ is a single line segment,
  for $r=\frac 2 {\sqrt 3} \rho$,
  and hence $S_1\cap \partial \Delta$ is composed by a single point, 
  which we call $H_1\in Y$. 
  Let $S'_1 := \overline{S_1\setminus \Delta}$, and notice that
  $S_1'$ is obtained from $S_1$ by removing the segment $[T,H_1)$.
  Hence, from $S_1\in \M(A_1\cup \ENCLOSE{T})$, we get $S'_1\in \M(A_1\cup \ENCLOSE{H_1})$. 
  Since $H_1\in \Delta$, we have also $S'_1\in \St(A_1 \cup \Delta)$.
  If we define $\Sigma'_1 = \overline{\Sigma\setminus \Delta}$,
  by Theorem~\ref{th:main} we know that $\Sigma'_1 \in \M(A_1\cup Y)$ 
  and hence $\Sigma'_1 \in \M(A_1\cup \Delta)$ because the projection of $A_1$ onto the line 
  $Y$ is contained in $\Delta$. 
  Hence $\H^1(\Sigma'_1)\le \H^1(S'_1)$.
  Repeat the reasoning for $k=2,3$ to define $H_k$, $S'_k$ and $\Sigma'_k$
  as we did for $k=1$.
  We obtain
  \begin{align*}
    \H^1(S) &= \H^1(S\cap \Delta) + \sum_{k=1}^3 \H^1(S'_k),\\
    \H^1(\Sigma_u) &= \H^1(\Sigma_u\cap \Delta) + \sum_{k=1}^3 \H^1(\Sigma'_k),\\
    \H^1(\Sigma'_k) & \le \H^1(S'_k) \qquad \text{for $k=1,2,3$}.
  \end{align*}
  Now notice that $S\cap \Delta = [T,H_1]\cup [T,H_2]\cup [T,H_3]$
  and hence, by Lemma~\ref{lm:tripod}, $\H^1(S\cap \Delta) \ge \H^1(\Sigma_u\cap \Delta)$,
  since $\abs{H_1 - T}\ge \dist(T,Y)$, and analogously for $k=2,3$.
  Putting all together, we obtain $\H^1(S)\ge \H^1(\Sigma_u)$
  which means that $\Sigma_u\in \M(A_u)$, and actually $H^1(S)=\H^1(\Sigma_u)$.

  Since we have the equality $\H^1(S)=\H^1(\Sigma_u)$ the inequalities 
  $\H^1(S\cap \Delta)\le \H^1(\Sigma_u\cap \Delta)$ and
  $\H^1(S'_k)\ge \H^1(\Sigma'_k)$, $k=1,2,3$, must all be equalities.
  It means that $S'_k\in \M(A_k\cup \Delta)$ for $k=1,2,3$, and,
  by Theorem~\ref{th:main}, we might conclude that $S'_k\cup A_k = \Sigma'_k$.
  As a consequence we have also $S\cap \Delta = \Sigma_u\cap \Delta$ because 
  the tripod $S\cap \Delta$ is uniquely determined by the points of $S$ on 
  the boundary of $\Delta$.
\end{proof}

The Theorem below shows that, in fact, the constructed set 
$\Sigma_u$ is a universal Steiner tree, 
in the sense that it contains homeomorphic copies of all 
Steiner trees.

\begin{theorem}
\label{th:universality}%
Let $\Sigma_u$ be the universal tree as defined in the introduction
with such $\lambda>0$ that $\Sigma_u\in \M(A_u)$.
Then the following properties hold.
\begin{enumerate}
  \item[(i)]
Given any $S$ connected subset of $\Sigma_u$, 
letting $D=\partial_{\Sigma_u}(S\setminus A_u)$ (the topological boundary of $S\setminus A_u$ relative 
to $\Sigma_u$),
one has $S\in \M(D)$.
  \item[(ii)]
Given any abstract, possibly infinite, binary tree (i.e.\ a graph 
such that every pair of vertices has a unique finite path connecting them and all vertices 
have order at most 3) 
there exists a subset $S\subset \Sigma_u\setminus A_u$
such that $S$ is homeomorphic to the topological space given by the abstract tree 
(and hence, by the previous point, $S\in \M(D)$ is a Steiner tree for 
the set of its endpoints $D$).
\item[(iii)] 
  Let $X$ be any metric space with the Heine-Borel property
  (for example $X=\RR^d$), let $D\subset X$ be a compact set 
  and $S\in \M(D)$ with $\H^1(S)<+\infty$.
  Then every connected component of $S\setminus D$ 
  is homeomorphic to a subtree of $\Sigma_u$.
\end{enumerate}
\end{theorem}

\begin{remark}
  Notice that the set $\Sigma$ is also universal in the weaker 
  sense of containing homeomerphic copies of every \emph{finite}
  Steiner tree, but, for example, it does not contain a 
  homeomorphic copy of $\Sigma_u$.
  This weaker form of universality has been used, e.g., 
  in~\cite{BasCheRasTep24}.
  \end{remark}
\begin{proof}[Proof of Theorem~\ref{th:universality}]
Clearly $\H^1(S)\le \H^1(\Sigma_u)<+\infty$.
To prove (i) take any set $S'$ such that $S'\cup D$ is connected 
and $\H^1(S')<+\infty$. 
Consider 
\[
  \Sigma' 
      := (\Sigma_u\setminus S) \cup S' \cup D.
\]

\emph{Step 1.} We claim that $\Sigma'$ is connected.
After fixing $x_0\in D$ we are going to show that for any $x\in \Sigma'$ there exists 
a connected set containing $x$ and $x_0$ and contained in $\Sigma'\cup D$.
If $x\in S'\cup D$ then just take $S'\cup D$ and we are done.
Otherwise 
\[
  x\in \Sigma'\setminus (S'\cup D) 
    = (\Sigma_u\setminus S)\setminus(S'\cup D)
    = \Sigma_u \setminus (S\cup S' \cup D),
\]
in particular $x\not\in S$ and hence $x\not \in S\setminus A_u$.
By Lemma~\ref{lemma_da_fare} applied to the connected component $C$ of 
$\Sigma_u\setminus(S\setminus A_u)$ containing $x$,
since $\Sigma_u$ is locally connected, we know that $\bar C\cap \partial (S\setminus A_u)$
is non empty. Hence there is $x_1\in \bar C \cap \partial (S\setminus A_u)$ and $\bar C$
is connected. Therefore $\bar C \cup (S'\cup D)$ is a connected set containing 
both $x$ and $x_0$, showing the claim.

\emph{Step 2.} 
We claim that $\H^1(D)=0$. 
In fact the points of $D$ which are not elements of $A_u$ are contained 
in one of the countably many segments which are the arcs of the tree $\Sigma_u\setminus A_u$.
If one of these segments contains two points of $D$ then $S$ would be completely contained 
in the segment and hence these are the only two points of $D$.
Otherwise the number of points of $D\setminus A_u$ is countable and hence $\H^1(D\setminus A_u)=0$.
On the other hand $\H^1(D\cap A_u)\le \H^1(A_u)=0$ which is true when $\lambda < \frac 1 2$
i.e.\ when $\Sigma_u$ is a tree.

\emph{Step 3.}
It is easy to check that $A_u\subset \Sigma'$. 
Since $\Sigma'$ is connected, we hence have that $\Sigma'$ is in fact 
a competitor to the minimality of $\Sigma_u\in \M(A_u)$ and thus $\H^1(\Sigma')\le \H^1(\Sigma_u)$.
But since 
\[
  \H^1(\Sigma') \le \H^1(\Sigma_u)-\H^1(S) + \H^1(S') + \H^1(D) 
\]
we obtain $\H^1(S)\le \H^1(S')$ 
which proves $S\in \M(D)$, concluding the proof of (i).

\emph{Step 4.} For the proof of (ii) we take any abstract binary tree $T$ and take an arbitrary 
vertex $v$ of $T$.
We map $v$ to the point $(0,0)\in \Sigma_u$. 
Then we consider the $k\le 3$ edges of $T$ meeting in $v$ and map them 
to $k$ of the three segments of $\Sigma_u$ meeting in $(0,0)$. 
We then proceed by induction on the graph distance of the vertices from $v$: once 
we have mapped all the vertices up to distance $n$, for each of these vertices 
we have at most two adjacent vertices at distance $n+1$ which we can map to one or both of the two 
adjacent vertices of the tree in $\Sigma_u$.
Since both $T$ and $\Sigma_u\setminus A_u$ are trees, 
the additional vertices discovered in the induction process were not already used.
And since $T$ is connected, all its vertices are reached by the induction.

\emph{Step 5.} For the point (iii), given any compact set $D\subset \RR^d$ and any $S\in \M(D)$,
we recall that (by Proposition~\ref{prop:PaoSte}) every connected component $S_0$ 
of $S\setminus D$ is a locally finite binary tree. 
Hence, by (ii), $S_0$ is homeomorphic to a subset of $\Sigma_u$.
\end{proof}

\begin{lemma}\label{lemma_da_fare}
Let $X$ be a locally connected metric space, $A\subset X$, and 
$C$ a connected component of $A$. Then $\partial C \subset \partial A$.
\end{lemma}
\begin{proof}
Clearly $\partial C\subset \bar C\subset \bar A$.
Take any $x\in \partial C$ and $U$ any neighbourhood of $x$ in $X$.
Let $V\subset U$ be a connected neighbourhood of $x$.
Since $C\cup\ENCLOSE{x}$ is connected 
then $C\cup V$ is also connected. But $C$ is a maximal connected subset of $A$ 
and therefore we conclude that $V\setminus A$ is non empty, hence $U\setminus A$ is non empty 
which means that $x\in \partial A$.
\end{proof}

\begin{remark}
Notice that if $X$ is not locally connected, Lemma~\ref{lemma_da_fare} is not true.
In fact if $X$ is the so called Knaster-Kuratowski fan, which is even connected, 
there exist $x\in X$ such that $A:= X\setminus\ENCLOSE{x}$ is totally disconnected hence
every connected component $C$ of $A$ is a singleton $C=\ENCLOSE{c}$.
Hence $\partial C = \ENCLOSE{c}$ while $\partial A = \ENCLOSE{x}$ and clearly $c\neq x$.
\end{remark}

\bibliographystyle{plain}
  
\end{document}